\documentclass[11pt]{amsart}
\usepackage{mathrsfs}
\usepackage{amssymb,latexsym}
\usepackage{pstricks,color}

 \numberwithin{equation}{section}

\newtheorem{theorem}{Theorem}[section]

\newtheorem{corollary}[theorem]{Corollary}

\newtheorem{remark}[theorem]{Remark}

\newtheorem{lemma}[theorem]{Lemma}

\def\qed{\hfill $\Box$}
\def\pf{\noindent {\it Proof.} }

\title{ Some properties of a class of refined Eulerian polynomials }

\begin{document}
\maketitle
\begin{center}
Yidong Sun$^{\dag}$\footnote{Corresponding author: Yidong Sun.} and Liting Zhai$^{\ddag}$

School of Science, Dalian Maritime University, 116026 Dalian, P.R. China\\[5pt]

{\it Email: $^{\dag}$sydmath@dlmu.edu.cn\ \ \   $^{\ddag}$zlt166@163.com }

\end{center}\vskip0.2cm

\subsection*{Abstract} In recent, H. Sun defined a new kind of refined Eulerian polynomials, namely,
\begin{eqnarray*}
A_n(p,q)=\sum_{\pi\in \mathfrak{S}_n}p^{{\rm odes}(\pi)}q^{{\rm edes}(\pi)}
\end{eqnarray*}
for $n\geq 1$, where ${odes}(\pi)$ and ${edes}(\pi)$ enumerate the number of descents of permutation $\pi$ in odd and even positions, respectively.
In this paper, we build an exponential generating function for $A_{n}(p,q)$ and establish an explicit formula for $A_{n}(p,q)$ in terms of Eulerian polynomials $A_{n}(q)$ and
$C(q)$, the generating function for Catalan numbers. In certain special case, we set up a connection between $A_{n}(p,q)$ and $A_{n}(p,0)$ or $A_{n}(0,q)$,
and express the coefficients of $A_{n}(0,q)$ by Eulerian numbers. Specially, this connection creates a new relation between Euler numbers and Eulerian numbers.

\medskip

{\bf Keywords}: Eulerian polynomial; Eulerian number; Euler number; Descent; Alternating permutation; Catalan number.

\noindent {\sc 2010 Mathematics Subject Classification}: Primary 05A05; Secondary 05A15, 05A10.

{\bf \section{ Introduction } }

Let $\mathfrak{S}_n$ denote the set of all permutations on $[n]=\{1,2,\dots,n\}$. For a permutation
$\pi= a_1a_2\dots a_n \in \mathfrak{S}_n$, an index $i \in [n-1]$ is a {\it descent} of $\pi$ if $a_i > a_{i+1}$, and $des(\pi)$ denotes
the number of descents of $\pi$. It is well known that $A_{n,k}$, the {\it Eulerian number} \cite[A008292]{Sloane},
counts the number of permutations $\pi\in \mathfrak{S}_n$ with $k$ descents and obeys the following recurrence \cite{Comtet}
\begin{eqnarray*}
A_{n,k} \hskip-.22cm &=& \hskip-.22cm (n-k)A_{n-1,k-1}+(k+1)A_{n-1,k}, \ \ (n>k\geq 0)
\end{eqnarray*}
with $A_{n,0}=1$ for $n\geq 0$ and $A_{n,k}=0$ for $1\leq n\leq k$ or $k<0$. The exponential generating function \cite{Comtet} for $A_{n,k}$ is
\begin{eqnarray}
\mathcal{E}(q;t)=1+\sum_{n\geq 1}A_n(q)\frac{t^n}{n!}= 1+\sum_{n\geq 1}\sum_{k=0}^{n-1}A_{n,k}q^{k}\frac{t^n}{n!}=\frac{1-q}{e^{t(q-1)}-q}.
\end{eqnarray}
where
\begin{eqnarray*}
A_n(q)=\sum_{\pi\in \mathfrak{S}_n}q^{{\rm des}(\pi)}=\sum_{k=0}^{n-1}A_{n,k}q^{k},
\end{eqnarray*}
is the classical {\it Eulerian polynomial} \cite{FoatSch}. The Eulerian polynomials have a rich history and appear 
in a large number of contexts in combinatorics; see \cite{Petersen} for a detailed exposition.

Recently, Sun \cite{HSun} introduced a new kind of {\it refined Eulerian polynomials} defined by
\begin{eqnarray*}
A_n(p,q)=\sum_{\pi\in \mathfrak{S}_n}p^{{\rm odes}(\pi)}q^{{\rm edes}(\pi)}
\end{eqnarray*}
for $n\geq 1$, where ${odes}(\pi)$ and ${edes}(\pi)$ enumerate the number of descents of permutation $\pi$ in odd and even positions, respectively.
The polynomial $A_n(p,q)$ is a bivariate polynomial of degree $n-1$, and the monomial with degree $n-1$ is exactly
$p^{\lfloor\frac{n}{2}\rfloor}q^{\lfloor\frac{n-1}{2}\rfloor}$,
where $\lfloor x\rfloor$ denotes the largest integer $\leq x$. When $p=q$, $A_n(p,q)$ reduces to the Eulerian polynomial $A_n(q)$.

For convenience, denote by
\begin{eqnarray*}
\tilde{A}_n(p,q)=\left\{
\begin{array}{ll}
A_n(p,q),      &  {\rm if }\ n=2m+1, \\[5pt]
(1+q)A_n(p,q), &  {\rm if }\ n=2m+2,
\end{array}
\right.
\end{eqnarray*}
for $n\geq 1$ and $\tilde{A}_0(p,q)=A_0(p,q)=1$. Sun \cite{HSun} showed that the
(modified) refined Eulerian polynomials $\tilde{A}_{n}(p,q)$ is
palindromic (symmetric) of draga ${\lfloor\frac{n}{2}\rfloor}$. She
also provided certain explicit formulas for special cases, namely,
\begin{eqnarray}\label{eqn 1.5}
A_n(p,1) \hskip-.22cm &=& \hskip-.22cm \frac{n!}{2^{\lfloor\frac{n}{2}\rfloor}}(1+p)^{\lfloor\frac{n}{2}\rfloor}, \label{eqn 1.2}\\
A_n(1,q) \hskip-.22cm &=& \hskip-.22cm \frac{n!}{2^{\lfloor\frac{n-1}{2}\rfloor}}(1+q)^{\lfloor\frac{n-1}{2}\rfloor}. \label{eqn 1.3}
\end{eqnarray}

Note that a permutation $\pi\in \mathfrak{S}_n$ such that ${\rm odes}(\pi)={\lfloor\frac{n}{2}\rfloor}$ and ${\rm edes}(\pi)=0$
(or ${\rm odes}(\pi)=0$ and ${\rm edes}(\pi)={\lfloor\frac{n-1}{2}\rfloor}$ ) is exactly an alternating (or reverse alternating) permutation.
Recall that a permutation $\pi=a_1a_2\cdots a_n \in \mathfrak{S}_n$ is {\it alternating} (or {\it reverse alternating} ) \cite{StanleyAl} if $a_1>a_2<a_3>\cdots $
(or $a_1<a_2>a_3<\cdots $). It is well known
that the {\it Euler number} $E_n$ \cite[A000111]{Sloane} counts the (reverse) alternating permutations in $\mathfrak{S}_n$,
which has the remarkable generating function \cite{StanleyAl},
\begin{eqnarray*}
\sum_{n\geq 0}E_n\frac{t^n}{n!} \hskip-.22cm &=& \hskip-.22cm \tan(t)+\sec(t) \\
                                \hskip-.22cm &=& \hskip-.22cm 1+ t+ \frac{t^2}{2!}+2\frac{t^3}{3!}+5\frac{t^4}{4!}+16\frac{t^5}{5!}+61\frac{t^6}{6!}+272\frac{t^7}{7!}+1385\frac{t^8}{8!}+\cdots.
\end{eqnarray*}
It produces that
\begin{eqnarray*}
\sum_{n\geq 0}E_{2n}\frac{t^{2n}}{(2n)!}       \hskip-.22cm &=& \hskip-.22cm \sec(t), \\
\sum_{n\geq 0}E_{2n+1}\frac{t^{2n+1}}{(2n+1)!} \hskip-.22cm &=& \hskip-.22cm \tan(t).
\end{eqnarray*}
For this reason $E_{2n}$ is sometimes called a {\it secant number} and $E_{2n+1}$ a {\it tangent number}.
See \cite{StanleyAl} for a survey of alternating permutations. The following associated generating functions are useful,
\begin{eqnarray*}
\sum_{n\geq 0}(-1)^{n}E_{2n}\frac{t^{2n}}{(2n)!}       \hskip-.22cm &=& \hskip-.22cm \frac{2}{e^t+e^{-t}}=\frac{2e^{t}}{e^{2t}+1}=\mathcal{E}(-1;t)e^{-t}, \\
\sum_{n\geq 0}(-1)^{n}E_{2n+1}\frac{t^{2n+1}}{(2n+1)!} \hskip-.22cm &=& \hskip-.22cm \frac{e^t-e^{-t}}{e^t+e^{-t}}=\frac{e^{2t}-1}{e^{2t}+1}=\mathcal{E}(-1;t)-1.
\end{eqnarray*}
It establishes a connection between Eulerian numbers and Euler numbers, that is,
\begin{eqnarray}
E_{2n+1}  \hskip-.22cm &=& \hskip-.22cm (-1)^{n}A_{2n+1}(-1), \label{eqn A1} \\
E_{2n+3}  \hskip-.22cm &=& \hskip-.22cm (-1)^{n}2A_{2n+2}'(-1), \label{eqn A2}
\end{eqnarray}
where $A_{n}'(q)$ is the derivative of $A_n(q)$ with respect to $q$.

In this paper, we further study the refined Eulerian polynomials $A_{n}(p,q)$ or $\tilde{A}_{n}(p,q)$.
The remainder of this paper is organized as follows. The next section will be devoted to
building an exponential generating function for $A_{n}(p,q)$ and to establishing an explicit
formula for $A_{n}(p,q)$ in terms of Eulerian polynomials $A_{n}(q)$ and
$C(q)=\frac{1-\sqrt{1-4q}}{2q}$, the generating function for Catalan numbers. The third section will discuss the special case $p=0$ or $q=0$, set up a connection between $A_{n}(p,q)$ and $A_{n}(p,0)$ or $A_{n}(0,q)$,
and express the coefficients of $A_{n}(0,q)$ by Eulerian numbers. Specially, it creates a new relation between Euler numbers and Eulerian numbers.

\vskip0.5cm
\section{The explicit formula for $A_n(p,q)$}

In this section, we consider the bivariate polynomials $A_n(p,q)$ and find an explicit formula for $A_n(p,q)$.
First, we need the following lemma.

\begin{lemma} For any integer $n\geq 1$, there holds
\begin{eqnarray}
A_{2n}(p,q)    \hskip-.22cm &=& \hskip-.22cm (1+p)A_{2n-1}(p,q)+(p+q)\sum_{i=1}^{n-1}\binom{2n-1}{2i-1}A_{2i-1}(p,q)A_{2n-2i}(p,q),      \label{eqn 1.6.1} \\
A_{2n+1}(p,q)  \hskip-.22cm &=& \hskip-.22cm A_{2n}(p,q)+p\sum_{i=0}^{n-1}\binom{2n}{2i}A_{2i}(p,q)A_{2n-2i}(q,p)           \label{eqn 1.6.2} \\
               & &   \hskip2cm  +q\sum_{i=1}^{n}\binom{2n}{2i-1}A_{2i-1}(p,q)A_{2n-2i+1}(p,q). \nonumber
\end{eqnarray}
\end{lemma}

\pf For any $\pi=a_1a_2\dots a_{2n} \in \mathfrak{S}_{2n}$, according to $a_k=2n$, $\pi$ can be partitioned into
$\pi=\pi_1(2n)\pi_2$ with $\pi_1=a_1a_2\dots a_{k-1}$ and $\pi_2=a_{k+1}a_{k+2}\dots a_{2n}$. Let $S=\{a_1, a_2,\dots, a_{k-1}\}$, then
$S$ is a $(k-1)$-subset of $[2n-1]$ and $\pi_2$ is a certain permutation of $[2n-1]-S$. If $\pi_2$ is empty, that is $a_{2n}=2n$, then
$\pi_1\in \mathfrak{S}_{2n-1}$, which are totally counted by $A_{2n-1}(p,q)$ according to ${\rm odes}(\pi_1)$ and ${\rm edes}(\pi_1)$.
If $\pi_2$ is not empty, that is $1\leq k< 2n$, then
\begin{eqnarray*}
{\rm odes}(\pi)= {\rm odes}(\pi_1)+{\rm odes}(\pi_2), \  {\rm edes}(\pi)= {\rm edes}(\pi_1)+{\rm edes}(\pi_2)+1,  \ {\rm when}\ k\ {\rm even},  \\
{\rm odes}(\pi)= {\rm odes}(\pi_1)+{\rm edes}(\pi_2)+1,\  {\rm edes}(\pi)= {\rm edes}(\pi_1)+{\rm odes}(\pi_2),   \ {\rm when}\ k\ {\rm odd}.
\end{eqnarray*}
Therefore, there are $\binom{2n-1}{k-1}$ choices to choose $S\in
[2n-1]$, all permutations $\pi_1$ of $S$ are counted by
$A_{k-1}(p,q)$ and all permutations $\pi_2$ of $[2n-1]-S$ are
counted by $A_{2n-k}(p,q)$, so all permutations $\pi=\pi_1(2n)\pi_2$
are counted by $A_{2i-1}(p,q)qA_{2n-2i}(p,q)$ when $k=2i$ for $1\leq
i< n$, and counted by $A_{2n-2i}(p,q)pA_{2i-1}(q,p)$ when $k=2(n-i)+1$ for
$1\leq i\leq n$. Note that $A_{k}(p,q)=A_{k}(q,p)$ when $k$ is odd.
To summarize all these cases, we obtain (\ref{eqn 1.6.1}).

Similarly, one can prove (\ref{eqn 1.6.2}), the details are omitted. \qed\vskip0.2cm

Let $A^{(e)}(p,q;t)$ and $A^{(o)}(p,q;t)$ be the exponential generating functions for $A_{2n}(p,q)$
and $A_{2n+1}(p,q)$ respectively, i.e.,
\begin{eqnarray*}
A^{(e)}(p,q;t) \hskip-.22cm &=& \hskip-.22cm \sum_{n\geq 1}A_{2n}(p,q)\frac{t^{2n}}{(2n)!}, \\
A^{(o)}(p,q;t) \hskip-.22cm &=& \hskip-.22cm \sum_{n\geq 0}A_{2n+1}(p,q)\frac{t^{2n+1}}{(2n+1)!}.
\end{eqnarray*}
 Then Lemma 2.1 suggests that
\begin{eqnarray}
\frac{\partial A^{(e)}(p,q;t)}{\partial t} \hskip-.22cm &=& \hskip-.22cm A^{(o)}(p,q;t)\big((1+p)+(p+q)A^{(e)}(p,q;t)\big), \label{eqnAE1}\\
\frac{\partial A^{(o)}(p,q;t)}{\partial t} \hskip-.22cm &=& \hskip-.22cm 1+ A^{(e)}(p,q;t)+ p(A^{(e)}(p,q;t)+1)A^{(e)}(q,p;t)+qA^{(o)}(p,q;t)^2. \label{eqnAO1}
\end{eqnarray}

Notice that $(1+q)A_{2n}(p,q)=(1+p)A_{2n}(q,p)$ and $A_{2n-1}(p,q)=A_{2n-1}(q,p)$ for $n\geq 1$,
one has $A^{(e)}(q,p;t)=\frac{1+q}{1+p}A^{(e)}(p,q;t)$ and $A^{(o)}(q,p;t)= A^{(o)}(p,q;t)$.
Then after simplification, (\ref{eqnAO1}) produces
{\small
\begin{eqnarray}
\frac{\partial A^{(o)}(p,q;t)}{\partial t} \hskip-.22cm &=& \hskip-.22cm \frac{1}{2}\frac{\partial (A^{(o)}(p,q;t)+A^{(o)}(q,p;t))}{\partial t} \nonumber  \\
                                           \hskip-.22cm &=& \hskip-.22cm 1+ (1+q)A^{(e)}(p,q;t)+ \frac{(1+q)(p+q)}{2(1+p)}A^{(e)}(p,q;t)^2+\frac{p+q}{2}A^{(o)}(p,q;t)^2. \label{eqnAO2} \hskip-1cm
\end{eqnarray} }

Let $y=\frac{p+q}{2(1+pq)}$ and $x=yC(y^2)$, where $C(y)=\frac{1-\sqrt{1-4y}}{2y}$ is 
the generating function of the Catalan numbers $C_n=\frac{1}{n+1}\binom{2n}{n}$ for $n\geq 0$.
By the relation $C(y)=1+yC(y)^2$, we have $y=\frac{x}{1+x^2}$, $\frac{2x}{1+x^2}=\frac{p+q}{1+pq}$, $\frac{2x}{(1+x)^2}=\frac{p+q}{(1+p)(1+q)}$ and
$\frac{(1+x)^2}{1+x^2}=\frac{(1+p)(1+q)}{1+pq}$. Define
\begin{eqnarray*}
B^{(e)}(p,q;t)  \hskip-.22cm &=& \hskip-.22cm \frac{1+x}{1+q}\left\{\frac{\mathcal{E}(x;\sqrt{\frac{1+pq}{1+x^2}}t)+\mathcal{E}(x;-\sqrt{\frac{1+pq}{1+x^2}}t)}{2}-1\right\}, \label{eqnBE1}\\
B^{(o)}(p,q;t)  \hskip-.22cm &=& \hskip-.22cm \frac{\mathcal{E}(x;\sqrt{\frac{1+pq}{1+x^2}}t)-\mathcal{E}(x;-\sqrt{\frac{1+pq}{1+x^2}}t)}{2\sqrt{\frac{1+pq}{1+x^2}}}. \label{eqnBO1}
\end{eqnarray*}
Now take partial derivative for $B^{(e)}(p,q;t)$ and $B^{(o)}(p,q;t)$ with respect to $t$, we obtain
\begin{eqnarray}
\frac{\partial B^{(e)}(p,q;t)}{\partial t} \hskip-.22cm &=& \hskip-.22cm \frac{(1+x)}{(1+q)}\frac{(1+pq)}{(1+x^2)}B^{(o)}(p,q;t)\Big((1+x)+\frac{2x(1+q)}{1+x}B^{(e)}(p,q;t)\Big), \nonumber \\
                                           \hskip-.22cm &=& \hskip-.22cm  B^{(o)}(p,q;t)\Big((1+p)+(p+q)B^{(e)}(p,q;t)\Big), \label{eqnBE2}
\end{eqnarray}
and
\begin{eqnarray}
\frac{\partial B^{(o)}(p,q;t)}{\partial t} \hskip-.22cm &=& \hskip-.22cm 1+(1+q)B^{(e)}(p,q;t) \nonumber \\
                                           & &  +\frac{x(1+q)^2}{(1+x)^2}B^{(e)}(p,q;t)^2+\frac{x(1+pq)}{1+x^2}B^{(o)}(p,q;t)^2,  \nonumber \\
                                           \hskip-.22cm &=& \hskip-.22cm 1+(1+q)B^{(e)}(p,q;t)+\frac{(1+q)(p+q)}{2(1+p)}B^{(e)}(p,q;t)^2+\frac{p+q}{2}B^{(o)}(p,q;t)^2. \label{eqnBO2}
\end{eqnarray}
By (\ref{eqnAE1}), (\ref{eqnAO1}) and (\ref{eqnBE2}), (\ref{eqnBO2}), one can see that $A^{(e)}(p,q;t)$, $A^{(o)}(p,q;t)$ and $B^{(e)}(p,q;t)$, $B^{(o)}(p,q;t)$ satisfy the
same differential equations of order one. On the other hand, it is routine to verify that for $1\leq n\leq 3$ the coefficients of $t^{2n}$ in $A^{(e)}(p,q;t)$ and $B^{(e)}(p,q;t)$ coincide,
as well as the ones of $t^{2n-1}$ in $A^{(o)}(p,q;t)$ and $B^{(o)}(p,q;t)$. Hence we obtain
the exponential generating functions of $A_{2n}(p,q)$ and $A_{2n-1}(p,q)$ for $n\geq 1$ as follows.
\begin{theorem} There hold
\begin{eqnarray}
A^{(e)}(p,q;t)  \hskip-.22cm &=& \hskip-.22cm \frac{1+x}{1+q}\left\{\frac{\mathcal{E}(x;\sqrt{\frac{1+pq}{1+x^2}}t)+\mathcal{E}(x;-\sqrt{\frac{1+pq}{1+x^2}}t)}{2}-1\right\}, \label{eqnAE3}\\
A^{(o)}(p,q;t)  \hskip-.22cm &=& \hskip-.22cm \frac{\mathcal{E}(x;\sqrt{\frac{1+pq}{1+x^2}}t)-\mathcal{E}(x;-\sqrt{\frac{1+pq}{1+x^2}}t)}{2\sqrt{\frac{1+pq}{1+x^2}}}. \label{eqnAO3}
\end{eqnarray}
where $x=yC(y^2)$, $y=\frac{p+q}{2(1+pq)}$ and $C(y)=\frac{1-\sqrt{1-4y}}{2y}$.
\end{theorem}

Comparing the coefficient of $\frac{t^n}{n!}$ in both sides of (\ref{eqnAE3}) and (\ref{eqnAO3}), we have the following explicit formula for $A_{n}(p,q)$.
\begin{corollary}\label{coro 2.3}
For any integer $n\geq 0$, there hold
{\small
\begin{eqnarray}
A_{2n+1}(p,q) \hskip-.22cm &=& \hskip-.22cm A_{2n+1}\Big(\frac{p+q}{2(1+pq)}C\Big(\frac{(p+q)^2}{4(1+pq)^2}\Big)\Big)\Big(\frac{1+pq}{C\Big(\frac{(p+q)^2}{4(1+pq)^2}\Big)}\Big)^{n},  \label{eqnABE4}   \\
A_{2n+2}(p,q) \hskip-.22cm &=& \hskip-.22cm \frac{1+\frac{p+q}{2(1+pq)}C\Big(\frac{(p+q)^2}{4(1+pq)^2}\Big)}{1+q}A_{2n+2}\Big(\frac{p+q}{2(1+pq)}C\Big(\frac{(p+q)^2}{4(1+pq)^2}\Big)\Big)\Big(\frac{1+pq}{C\Big(\frac{(p+q)^2}{4(1+pq)^2}\Big)}\Big)^{n+1},  \label{eqnACE4} \hskip-1cm
\end{eqnarray}}
where $C(y)=\frac{1-\sqrt{1-4y}}{2y}$.
\end{corollary}

In the case $q=1$ in (\ref{eqnABE4}) and (\ref{eqnACE4}), by $C(\frac{1}{4})=2$ and $A_{n}(1)=n!$, one has
\begin{eqnarray*}
A_{2n+1}(p,1) \hskip-.22cm &=& \hskip-.22cm  \frac{(2n+1)!}{2^{n}}(1+p)^{n},   \\
A_{2n+2}(p,1) \hskip-.22cm &=& \hskip-.22cm  \frac{(2n+2)!}{2^{n+1}}(1+p)^{n+1},
\end{eqnarray*}
which is equivalent to (\ref{eqn 1.2}).

Similarly, the case $p=1$ in (\ref{eqnABE4}) and (\ref{eqnACE4}) also generates an equivalent form to (\ref{eqn 1.3}).

\vskip0.5cm

\section{The special case $p=0$ or $q=0$ }

In this section, we concentrate on the special case $p=0$ or $q=0$.

Since the symmetry of $\tilde{A}_{n}(p,q)$ for $n\geq 1$, there has
\begin{eqnarray}\label{eqn 3.1}
A_n(p,0)=\left\{
\begin{array}{ll}
A_n(0,p),      &  {\rm if }\ n=2m+1, \\[5pt]
(1+p)A_n(0,p), &  {\rm if }\ n=2m+2.
\end{array}
\right.
\end{eqnarray}

In fact, one can represent $\tilde{A}_{n}(p,q)$ in terms of $A_{n}(p,0)$ or $A_{n}(0,q)$.

\begin{theorem}\label{theo 3.1}
For any integer $n\geq 1$, there holds
\begin{eqnarray}\label{eqn 3.2}
\tilde{A}_n(p,q) \hskip-.22cm &=& \hskip-.22cm (1+pq)^{\lfloor\frac{n}{2}\rfloor}A_n\Big(\frac{p+q}{1+pq},0\Big).
\end{eqnarray}
Equivalently, for $n\geq 0$ there have
\begin{eqnarray}
A_{2n+1}(p,q)   \hskip-.22cm &=& \hskip-.22cm (1+pq)^{n}A_{2n+1}\Big(0, \frac{p+q}{1+pq}\Big), \label{eqn3.3A} \\
A_{2n+2}(p,q)   \hskip-.22cm &=& \hskip-.22cm (1+p)(1+pq)^{n}A_{2n+2}\Big(0, \frac{p+q}{1+pq}\Big). \label{eqn3.3B}
\end{eqnarray}
\end{theorem}

\pf The case $q=0$ in Corollary 2.3 gives rise to
\begin{eqnarray}
A_{2n+1}(p,0) \hskip-.22cm &=& \hskip-.22cm A_{2n+1}\Big(\frac{p}{2}C\Big(\frac{p^2}{4}\Big)\Big)C\Big(\frac{p^2}{4}\Big)^{-n},  \label{eqn 3.3}  \\
A_{2n+2}(p,0) \hskip-.22cm &=& \hskip-.22cm \Big(1+\frac{p}{2}C\Big(\frac{p^2}{4}\Big)\Big)A_{2n+2}\Big(\frac{p}{2}C\Big(\frac{p^2}{4}\Big)\Big)C\Big(\frac{p^2}{4}\Big)^{-n-1}. \label{eqn 3.4}
\end{eqnarray}
Then resetting $p:=\frac{p+q}{1+pq}$, by Corollary 2.3, we get
\begin{eqnarray*}
A_{2n+1}\Big(\frac{p+q}{1+pq},0\Big) \hskip-.22cm &=& \hskip-.22cm A_{2n+1}(p,q)(1+pq)^{-n}=\tilde{A}_{2n+1}(p,q)(1+pq)^{-n},    \\
A_{2n+2}\Big(\frac{p+q}{1+pq},0\Big) \hskip-.22cm &=& \hskip-.22cm (1+q)A_{2n+2}(p,q)(1+pq)^{-n-1}=\tilde{A}_{2n+2}(p,q)(1+pq)^{-n-1},
\end{eqnarray*}
which is equivalent to (\ref{eqn 3.2}), or by (\ref{eqn 3.1}), equivalent to (\ref{eqn3.3A}) and (\ref{eqn3.3B}). \qed\vskip.2cm

\begin{remark}
As stated by Sun \cite{HSun}, the refined Eulerian polynomials $\tilde{A}_{n}(p,q)$ can be expanded in terms of gamma basis, that is,
\begin{eqnarray*}
\tilde{A}_n(p,q) \hskip-.22cm &=& \hskip-.22cm \sum_{j=0}^{\lfloor\frac{n}{2}\rfloor}c_{n,j}(p+q)^j(1+pq)^{\lfloor\frac{n}{2}\rfloor-j},
\end{eqnarray*}
and she conjectured that for any $n\geq 1$, all $c_{n,j}$ are positive integers. From Theorem \ref{theo 3.1},
$$A_{n}(p,0)=\sum_{j=0}^{\lfloor\frac{n}{2}\rfloor}c_{n,j}p^{j},$$
it is obvious that $c_{n,j}$ are positive because $c_{n,j}$ counts the number of
permutations $\pi\in \mathfrak{S}_n$ with $j$ odd descents and with no even descents, and one can easily construct such a permutation $\pi=\pi_1(2j+1)(2j+2)\cdots n$,
where $\pi_1$ is an alternating permutation on $[2j]$.

\end{remark}

Let $a_{n,k}$ be the number of permutations $\pi\in \mathfrak{S}_n$ with $k$ even descents and with no odd descents, then
$$A_{n}(0,q)=\sum_{k=0}^{\lfloor\frac{n-1}{2}\rfloor}a_{n,k}q^{k},\ \ (n\geq 1), $$
by (\ref{eqn 3.1}), it's clear that $c_{n,k}=a_{n,k}$ when $n$ odd and $c_{n,k}=a_{n,k}+a_{n,k-1}$ when $n$ even.
Now we can establish several connections between Eulerian numbers $A_{n,k}$ and $a_{n,k}$.

\begin{corollary} For any integers $n\geq k\geq 0$, there hold
\begin{eqnarray}
A_{2n+1,k}  \hskip-.22cm &=& \hskip-.22cm \sum_{i=0}^{\lfloor\frac{k}{2}\rfloor}\binom{n-k+2i}{i}2^{k-2i}a_{2n+1,k-2i}, \label{coro3eqna1}\\
A_{2n+2,k}  \hskip-.22cm &=& \hskip-.22cm \sum_{{2i+j=k\ {\rm or} \ k-1}\atop {i,j\geq 0}}\binom{n-j}{i}2^{j}a_{2n+2,j}.\label{coro3eqna2}
\end{eqnarray}
Equivalently,
\begin{eqnarray}
a_{2n+1,k}  \hskip-.22cm &=& \hskip-.22cm \frac{1}{2^{k}}\sum_{i=0}^{\lfloor\frac{k}{2}\rfloor}(-1)^{i}\frac{n-k+2i}{n-k+i}\binom{n-k+i}{i}A_{2n+1,k-2i}, \label{coro3eqna3}\\
a_{2n+2,k}  \hskip-.22cm &=& \hskip-.22cm \frac{1}{2^{k}}\sum_{{2i+j+r=k}\atop {i,j,r\geq 0}}(-1)^{r+i}\frac{n-k+2i}{n-k+i}\binom{n-k+i}{i}A_{2n+2,j}.\label{coro3eqna4}
\end{eqnarray}
Specially, Euler numbers can be represented by Eulerian numbers as follows
\begin{eqnarray}
E_{2n+1}  \hskip-.22cm &=& \hskip-.22cm  a_{2n+1,n} =\frac{1}{2^{n}}\Big(-A_{2n+1,n}+2\sum_{i=0}^{\lfloor\frac{n}{2}\rfloor}(-1)^{i}A_{2n+1,n-2i}\Big), \label{coro3eqna5}\\
E_{2n+2}  \hskip-.22cm &=& \hskip-.22cm  a_{2n+2,n} =\frac{1}{2^{n}}\Big(-\sum_{j=0}^{n}(-1)^{n-j}A_{2n+2,j}+2\sum_{{0\leq 2i+j\leq n}\atop {i,j\geq 0}}(-1)^{n-i-j}A_{2n+2,j}\Big).\label{coro3eqna6} \hskip-1cm
\end{eqnarray}

\end{corollary}

\pf By (\ref{eqn3.3A}) and (\ref{eqn3.3B}), the case $p=q$ produces
\begin{eqnarray}
A_{2n+1}(q) \hskip-.22cm &=& \hskip-.22cm  (1+q^2)^{n}A_{2n+1}\Big(0,\frac{2q}{1+q^2}\Big)
                                            = \sum_{j=0}^{n}a_{2n+1,j}(2q)^{j}(1+q^2)^{n-j} ,  \label{eqn AB1} \\
A_{2n+2}(q) \hskip-.22cm &=& \hskip-.22cm  (1+q)(1+q^2)^{n}A_{2n+2}\Big(0,\frac{2q}{1+q^2}\Big)
                                            = (1+q)\sum_{j=0}^{n}a_{2n+2,j}(2q)^{j}(1+q^2)^{n-j}.  \label{eqn AB2} \hskip-.5cm
\end{eqnarray}
Comparing the coefficients of $q^{k}$ in both sides of (\ref{eqn AB1}) and (\ref{eqn AB2}), we can get (\ref{coro3eqna1}) and (\ref{coro3eqna2}).

By (\ref{eqn 3.1}) and (\ref{eqn 3.4}), using the series expansion \cite{Stanley},
\begin{eqnarray*}
C(t)^{\alpha} \hskip-.22cm &=&\hskip-.22cm \sum_{i\geq 0}\frac{\alpha}{2i+\alpha}\binom{2i+\alpha}{i} t^i,
\end{eqnarray*}
we have
\begin{eqnarray}
A_{2n+1}(0,q) \hskip-.22cm &=& \hskip-.22cm A_{2n+1}(q,0)=A_{2n+1}\Big(\frac{q}{2}C\Big(\frac{q^2}{4}\Big)\Big)C\Big(\frac{q^2}{4}\Big)^{-n}   \nonumber\\
              \hskip-.22cm &=& \hskip-.22cm \sum_{j=0}^{2n}A_{2n+1,j}\Big(\frac{q}{2}\Big)^{j}C\Big(\frac{q^2}{4}\Big)^{j-n}      \nonumber \\
              \hskip-.22cm &=& \hskip-.22cm \sum_{j=0}^{2n}A_{2n+1,j}\sum_{i\geq 0}\frac{j-n}{2i+j-n}\binom{2i+j-n}{i}\Big(\frac{q}{2}\Big)^{2i+j}   \nonumber \\
              \hskip-.22cm &=& \hskip-.22cm \sum_{j=0}^{2n}\sum_{i\geq 0}A_{2n+1,j}(-1)^{i}\frac{n-j}{n-i-j}\binom{n-i-j}{i}\Big(\frac{q}{2}\Big)^{2i+j},    \label{eqn ABC1}
\end{eqnarray}
and
\begin{eqnarray}
A_{2n+2}(0,q) \hskip-.22cm &=& \hskip-.22cm \frac{A_{2n+2}(q,0)}{1+q}=\frac{1+\frac{q}{2}C\Big(\frac{q^2}{4}\Big)}{1+q}A_{2n+2}\Big(\frac{q}{2}C\Big(\frac{q^2}{4}\Big)\Big)C\Big(\frac{q^2}{4}\Big)^{-n-1}   \nonumber\\
              \hskip-.22cm &=& \hskip-.22cm \frac{1}{1+\frac{q}{2}C\Big(\frac{q^2}{4}\Big)}A_{2n+2}\Big(\frac{q}{2}C\Big(\frac{q^2}{4}\Big)\Big)C\Big(\frac{q^2}{4}\Big)^{-n}   \nonumber\\
              \hskip-.22cm &=& \hskip-.22cm \sum_{j=0}^{2n+1}A_{2n+2,j}\sum_{r\geq 0}(-1)^{r}\Big(\frac{q}{2}\Big)^{j+r}C\Big(\frac{q^2}{4}\Big)^{r+j-n}      \nonumber \\
              \hskip-.22cm &=& \hskip-.22cm \sum_{j=0}^{2n+1}A_{2n+2,j}\sum_{r\geq 0}(-1)^{r}\sum_{i\geq 0}\frac{r+j-n}{2i+j+r-n}\binom{2i+r+j-n}{i}\Big(\frac{q}{2}\Big)^{2i+j+r}   \nonumber \\
              \hskip-.22cm &=& \hskip-.22cm \sum_{j=0}^{2n+1}\sum_{r,i\geq 0}A_{2n+2,j}(-1)^{r+i}\frac{n-j-r}{n-i-j-r}\binom{n-i-j-r}{i}\Big(\frac{q}{2}\Big)^{2i+j+r}.   \label{eqn ABC2}
\end{eqnarray}
Then taking the coefficients of $t^{k}$ in both sides of (\ref{eqn ABC1}) and (\ref{eqn ABC2}), we get (\ref{coro3eqna3}) and (\ref{coro3eqna4}).

Note that $a_{n,{\lfloor\frac{n-1}{2}\rfloor}}$ is the number of permutations $\pi\in \mathfrak{S}_n$ with ${\lfloor\frac{n-1}{2}\rfloor}$ even descents and with no odd descents,
such $\pi$ are exactly the reverse alternating permutations and vice verse. Clearly, $a_{n,{\lfloor\frac{n-1}{2}\rfloor}}=E_n$ for $n\geq 1$.
Setting $k=n$ in (\ref{coro3eqna3}) and (\ref{coro3eqna4}), after simplification we obtain (\ref{coro3eqna5}) and (\ref{coro3eqna6}). \qed\vskip.2cm

\begin{lemma}\label{lemmaABC}
For any integer $n\geq 1$, there holds
\begin{eqnarray}\label{eqn ABCD1}
A_n(0,-1)=\left\{
\begin{array}{ll}
(-1)^{m}\frac{E_{n}}{2^{m}},         &  {\rm if }\ n=2m+1, \\[5pt]
(-1)^{m}\frac{E_{n+1}}{2^{m+1}},     &  {\rm if }\ n=2m+2.
\end{array}
\right.
\end{eqnarray}
\end{lemma}
\pf Setting $p=q=-1$ in (\ref{eqn3.3A}) and (\ref{eqn3.3B}), by (\ref{eqn A1}) and (\ref{eqn A2}) we have
\begin{eqnarray*}
A_{2n+1}(0,-1) \hskip-.22cm &=& \hskip-.22cm  2^{-n}A_{2n+1}(-1)=(-1)^{n}\frac{E_{2n+1}}{2^n}, \\
A_{2n+2}(0,-1) \hskip-.22cm &=& \hskip-.22cm  2^{-n}\lim_{q\rightarrow -1} \frac{A_{2n+2}(q)}{1+q}
                  =2^{-n}A_{2n+2}'(-1)=(-1)^{n}\frac{E_{2n+3}}{2^{n+1}},
\end{eqnarray*}
which is equivalent to (\ref{eqn ABCD1}). \qed\vskip.2cm

The case $q=-1$ in (\ref{eqn3.3A}) and (\ref{eqn3.3B}), together with Lemma \ref{lemmaABC}, leads to
\begin{corollary}\label{coro3.5}
For any integer $n\geq 0$, there holds
\begin{eqnarray*}
A_{2n+1}(p,-1) \hskip-.22cm &=& \hskip-.22cm  \Big(\frac{p-1}{2}\Big)^{n}E_{2n+1},     \\
A_{2n+2}(p,-1) \hskip-.22cm &=& \hskip-.22cm  \frac{p+1}{2}\Big(\frac{p-1}{2}\Big)^{n}E_{2n+3}.
\end{eqnarray*}
\end{corollary}
The case $p=3$ in Corollary \ref{coro3.5} produces two settings counted by tangent numbers $E_{2n+1}$. See \cite{StanleyAl} for
further information on various combinatorial interpretations of Euler numbers $E_{n}$.

Let $\pi=a_1a_2\cdots a_{n-1}a_{n}\in \mathfrak{S}_n$, define the {\it reversal} of $\pi$ to be
$\pi^{r}=a_{n}a_{n-1}\cdots a_2a_1$, the {\it complement} of $\pi$ to be $\pi^{c}=(n+1-a_1)(n+1-a_2)\cdots (n+1-a_{n-1})(n+1-a_{n})$
and the {\it reversal-complement} of $\pi$ to be $\pi^{rc}:=(\pi^{r})^{c}=(\pi^{c})^{r}$. For any $\sigma\in \mathfrak{S}_{2n+1}$, note that
$i$ is a descent of $\sigma$ if and only if $2n+1-i$ is a descent of $\sigma^{rc}$. Specially, $i$ is an odd (even) descent of $\sigma$ if and only
if $2n+1-i$ is an even (odd) descent of $\sigma^{rc}$. Now we can return to consider the recurrence relations for $a_{n,k}$.

\begin{theorem} For any integers $n\geq k\geq 1$, there hold
\begin{eqnarray}
a_{2n,k}    \hskip-.22cm &=& \hskip-.22cm (n-k)a_{2n-1,k-1}+(k+1)a_{2n-1,k},   \label{eqn RecA} \\[5pt]
a_{2n+1,k}  \hskip-.22cm &=& \hskip-.22cm (n-k+1)a_{2n,k-2}+(n+1)a_{2n,k-1}+(k+1)a_{2n,k},  \label{eqn RecB}
\end{eqnarray}
with $a_{n,0}=1$ and $a_{n,\lfloor\frac{n-1}{2}\rfloor}=E_n.$
\end{theorem}

\pf Let $\alpha_{n,k}$ be the set of permutations $\pi\in \mathfrak{S}_n$ with $k$ even descents and with no odd descents,
so $|\alpha_{n,k}|=a_{n,k}$. In the $k=0$ case, $\alpha_{n,0}$ only contains the natural permutation $123\cdots (2n-1)(2n)$, and
in the $k=\lfloor\frac{n-1}{2}\rfloor$ case $\alpha_{n,\lfloor\frac{n-1}{2}\rfloor}$ is exactly the set of reverse alternating permutations
in $\mathfrak{S}_n$. This implies $a_{n,0}=1$ and $a_{n,\lfloor\frac{n-1}{2}\rfloor}=E_n.$

Set $X=\{2i|1\leq i\leq n-1\}$ and denote by $Des(\pi)$ the descent set of $\pi\in \mathfrak{S}_n$.
For any $\pi\in \alpha_{2n,k}$, $Des(\pi)$ is a $k$-subset of $X$. In order to prove (\ref{eqn RecA}),
there are exactly three cases should be considered, i.e.,

{\bf Case 1:} Given $\pi=a_1a_2\cdots a_{2n-2}a_{2n-1}\in \alpha_{2n-1,k}$, let
$\pi_1^{*}=\pi a_{2n}$ with $a_{2n}=2n$, one can easily check that $\pi_1^{*}\in \alpha_{2n,k}$ and $\pi$ can be
recovered by deleting $a_{2n}=2n$ in $\pi_1^{*}$. So there are totally $a_{2n-1,k}$
number of such $\pi_1^{*}$'s in the set $\alpha_{2n,k}$.

{\bf Case 2:} Given $\pi=a_1a_2\cdots a_{2j}a_{2j+1}\cdots a_{2n-2}a_{2n-1}\in \alpha_{2n-1,k}$ with $2j\in Des(\pi)$, define
$\pi_2^{*}=a_{2j+1}\cdots a_{2n-2}a_{2n-1}(2n)a_1a_2\cdots a_{2j}$. Clearly, subject to
$a_{2j}>a_{2j+1}$, we obtain $\pi_2^{*}\in \alpha_{2n,k}$ and vice versa. 
In these cases, there are totally $ka_{2n-1,k}$ contributions to the set $\alpha_{2n,k}$.

{\bf Case 3:} Given $\pi=a_1a_2\cdots a_{2j}a_{2j+1}\cdots a_{2n-2}a_{2n-1}\in \alpha_{2n-1,k-1}$ with $2j\in X-Des(\pi)$, define
$\pi_3^{*}=a_{2j+1}\cdots a_{2n-2}a_{2n-1}(2n)a_1a_2\cdots a_{2j}$. Similarly, subject to
$a_{2j}<a_{2j+1}$, we have $\pi_3^{*}\in \alpha_{2n,k}$ and vice versa.
In this case, there are totally $(n-k)a_{2n-1,k-1}$ contributions to the set $\alpha_{2n,k}$.

Hence, summarizing the above three cases generates (\ref{eqn RecA}) immediately.

In order to prove (\ref{eqn RecB}), by the relation $c_{n,k}=a_{n,k}+a_{n,k-1}$ when $n$ even, we need take the equivalent form into account,
\begin{eqnarray}
a_{2n+1,k}  \hskip-.22cm &=& \hskip-.22cm (n-k+1)c_{2n,k-1}+kc_{2n,k}+a_{2n,k}.  \label{eqn RecC}
\end{eqnarray}

Let $\beta_{n,k}$ be the set of permutations $\theta\in \mathfrak{S}_n$ with $k$ odd descents and with no even descents,
so $|\beta_{n,k}|=c_{n,k}$. Set $Y=\{2i-1|1\leq i\leq n\}$, for any $\theta\in \beta_{2n,k}$, $Des(\theta)$ is a $k$-subset of $Y$.
Similarly, there are precisely three cases to be taken consideration.

{\bf Case I:} Given $\theta=b_1b_2\cdots b_{2n}\in \alpha_{2n,k}$, let
$\theta_1^{*}=\theta b_{2n+1}$ with $b_{2n+1}=2n+1$, obviously $\theta_1^{*}\in \alpha_{2n+1,k}$ and $\theta$ can be
easily obtained by removing $b_{2n+1}=2n+1$ in $\theta_1^{*}$. Then there are exactly $a_{2n,k}$ number of such $\theta_1^{*}$'s in the set $\alpha_{2n+1,k}$.

{\bf Case II:} Given $\theta=b_1b_2\cdots b_{2j-1}b_{2j}b_{2j+1}\cdots b_{2n}\in \beta_{2n,k}$ with $2j-1\in Des(\theta)$, define
$\theta_2^{*}=(b_1b_2\cdots b_{2j-1})^{rc}(2n+1)b_{2j}'b_{2j+1}'\cdots b_{2n}'=b_1'b_2'\cdots b_{2j-1}'(2n+1)b_{2j}'b_{2j+1}'\cdots b_{2n}'$,
where $b_{i}'=2n+1-b_{2j-i}$ for $1\leq i\leq 2j-1$ and $b_{2j}'b_{2j+1}'\cdots b_{2n}'$ is a permutation on $[2n]-\{b_1', b_2', \dots, b_{2j-1}'\}$
which has the same relative order as $b_{2j}b_{2j+1}\cdots b_{2n}$. This procedure is naturally invertible subject to
$b_{2j-1}>b_{2j}$. Clearly, $2\ell-1$ with $1\leq \ell<j$ is an odd descent of $\theta$
if and only if $2j-2\ell$ is an even descent of $\theta_2^{*}$, and $2\ell-1$ with $j\leq \ell\leq n$
is an odd descent of $\theta$ if and only if $2\ell$ is an even descent of $\theta_2^{*}$. Then we get
$\theta_2^{*}\in \alpha_{2n+1,k}$. In this case, there are totally $kc_{2n,k}$ contributions to the set $\alpha_{2n+1,k}$.
For example, let $\theta=4\ 1\ 8\ 3\ 5\ 7\ 10\ 2\ 6\ 9\in \beta_{10,3}$ and $Des(\theta)=\{1, 3, 7\}$, we have three $\theta_2^{*}\in \alpha_{11,3}$, namely,
$$7\ 11\ 1\ 8\ 3\ 4\ 6\ 10\ 2\ 5\ 9,\hskip.3cm  3\ 10\ 7\ 11\ 2\ 4\ 6\ 9\ 1\ 5\ 8,\hskip.3cm  1\ 4\ 6\ 8\ 3\ 10\ 7\ 11\ 2\ 5\ 9.$$

{\bf Case III:} Given $\theta=b_1b_2\cdots b_{2j-1}b_{2j}b_{2j+1}\cdots b_{2n}\in \beta_{2n,k-1}$ with $2j-1\in Y-Des(\theta)$, define
$\theta_3^{*}=(b_1b_2\cdots b_{2j-1})^{rc}(2n+1)b_{2j}'b_{2j+1}'\cdots b_{2n}'=b_1'b_2'\cdots b_{2j-1}'(2n+1)b_{2j}'b_{2j+1}'\cdots b_{2n}'$,
where $b_{i}'=2n+1-b_{2j-i}$ for $1\leq i\leq 2j-1$ and $b_{2j}'b_{2j+1}'\cdots b_{2n}'$ is a permutation on $[2n]-\{b_1', b_2', \dots, b_{2j-1}'\}$
which has the same relative order as $b_{2j}b_{2j+1}\cdots b_{2n}$. This procedure is also invertible subject to
$b_{2j-1}<b_{2j}$. Similar to Case II, $\theta_3^{*}\in \alpha_{2n+1,k}$. In this case, there are totally $(n-k+1)c_{2n,k-1}$
contributions to the set $\alpha_{2n+1,k}$.

Hence, summing over all the three cases yields (\ref{eqn RecC}) immediately. \qed

Table 1.1 illustrates this triangle for $n$ up to $10$ and $k$ up to $4$.
\begin{center}
\begin{eqnarray*}
\begin{array}{|c|ccccc|}\hline
n/k & 0   & 1    & 2    & 3     & 4     \\\hline
  1 & 1   &      &      &       &         \\
  2 & 1   &      &      &       &          \\
  3 & 1   & 2    &      &       &          \\
  4 & 1   & 5    &      &       &         \\
  5 & 1   & 13   & 16   &       &        \\
  6 & 1   & 28   & 61   &       &        \\
  7 & 1   & 60   & 297  & 272   &        \\
  8 & 1   & 123  & 1011 & 1385  &        \\
  9 & 1   & 251  & 3651 & 10841 & 7936   \\
  10& 1   & 506  & 11706& 50666 & 50521  \\\hline
\end{array}
\end{eqnarray*}
Table 1.1. The first values of $a_{n,k}$.
\end{center}

\vskip.2cm

\section*{Acknowledgements} {The authors are grateful to the referees for
the helpful suggestions and comments. The Project is sponsored by ``Liaoning
BaiQianWan Talents Program".}

\vskip.2cm


\end{document}